%

\input amstex
\documentstyle{amsppt}
\magnification=\magstep1 
\pageheight{19cm}
 
\def\force{\Vdash}  
\def\notforce{\nVdash}
\define\k{\kappa}
\define\a{\alpha}
\redefine\l{\lambda}
\define\tr{\vartriangleleft}
\define\tS{\tr_S}

\define\ra{\rightarrow}

\define\res{\restriction}
\define\Uarr{\overrightarrow{\Cal U}}
\define\oU{o^{\Cal U}}
\define\lqot{\text{``}}
\define\rqot{\text{''}}
\define\Ult{\text{{\rm Ult}}}
\define\Gen{\text{{\rm Gen}}}

\redefine\crit{\text{{\rm crit}}\,}
\define\Reg{\text{{\rm Reg}}}
\define\Sing{\text{{\rm Sing}}}
\define\Club{\text{{\rm Club}}}
\define\NS{\text{{\rm NS}}}
\define\Tr{\text{{\rm Tr}}\,}

\topmatter
\title
Possible Behaviours of the Reflection \\
Ordering of Stationary Sets
\endtitle
\author  Ji\v r\' i Witzany   \endauthor
\thanks I want to express my gratitude to T.Jech
and also to J.Zapletal for many valuable discussions and 
remarks on the subject.
   \endthanks
\affil     The Pennsylvania State University  and
Charles University (Prague)                   \endaffil
\address Department of Mathematics, The Pennsylvania State University,
University Park, PA 16802                       \endaddress

\abstract
If $S,T$ are stationary subsets of a regular uncountable cardinal $\kappa$,
we say that $S$ reflects fully in $T$, $S<T$, if for almost all
$\alpha \in T$ (except a nonstationary set) $S \cap \alpha$ 
is stationary in $\alpha .$
This relation is known to be a well-founded partial ordering.
We say that a given poset $P$ is realized by the reflection ordering if there
is a maximal antichain $\langle X_p ; p \in P \rangle$ of stationary
subsets of $\Reg(\kappa)$ so that
$$\forall p,q \in P \; \forall S\subseteq X_p, T\subseteq X_q 
  \text{ stationary}:(S<T  \leftrightarrow p<_P q ) .$$

We prove that if $V=L[\Uarr],$ $\oU(\k)=\k^{++},$ and $P$ is an arbitrary 
well-founded poset of cardinality $\leq \k^+$
then there is a generic extension where P is realized by the reflection
ordering on $\kappa .$

\endabstract 
\date February 3, 1994 \enddate
\email witzany\@math.psu.edu                    \endemail 
\endtopmatter

\def\lheadline{\folio\hfil {\eightpoint JI\v R\' I WITZANY}\hfil}
\def\rheadline{\hfil{\eightpoint REFLECTION ORDERING OF STATIONARY SETS} \hfil\folio}
\headline={\ifodd\pageno\rheadline \else\lheadline\fi}

\document

\subhead
1. Introduction
\endsubhead
\vskip8pt

If $S$ is a stationary subset of a regular uncountable cardinal $\k$ then
{\it the trace of} $S$ is the set

$$ \Tr (S) = \{ \a < \k ; \; \; S \cap \a \; \text{is stationary
in}  \; \a \} $$
and we say that $S$ {\it reflects at} $\a\in \Tr(S)$.  If $S$ and $T$ are
both stationary,
we define
$$ S < T \;\; \text{if for almost all} \; \a \in T, \; \; \a \in
\Tr(S) $$
and say that $S$ {\it reflects fully} in $T$.
(Throughout the paper, ``for almost all" means ``except for a
nonstationary set of
points"). It can be proved that this relation is a well-founded
partial ordering (see [JW93] or [J84]).
 {\it The order} $o(S)$ of a stationary set of regular cardinals is defined
as the rank of S in the relation $<$:
$$ o(S)=\sup \{ o(T)+1 ; \; T \subseteq \Reg(\k ) \text{ is stationary and }  
T<S \} .$$
For a stationary set $T$ such that $T\cap \Sing(\k )$ is stationary 
define $o(T)=-1$. {\it The order of } $\k$ is then defined as
$$ o(\k) = \sup\{ o(S)+1 \; ;  S \subseteq \k \; \text{is stationary} \}.$$
Note that if $\Tr(S)$, where $S \subseteq \Reg(\k)$, is stationary then
$o(S) < o(\Tr(S))$ as $S<\Tr(S)$. It follows from [J84] that the order
$o(\k)$ provides a natural generalization of the Mahlo hierarchy: $\k$ is
exactly $o(\k)$-Mahlo if $o(\k) < \k^{+}$ and  greatly Mahlo
if $o(\k) \geq \k^{+}$.

Let $P$ be a well-founded poset, we say that {\it the reflection order} $<$
{\it realizes} $P$ if there is a maximal antichain $\langle
X_p ; p\in P \rangle $ of stationary subsets of $\Reg(\k)$ so that 
$$\forall p,q \in P \; \forall S\subseteq X_p, T\subseteq X_q 
\text{ stationary}:(S<T 
  \leftrightarrow p<_P q ) .$$
If $|P|\leq \k$ then it follows that
for any stationary $S\subseteq \Reg(\k)$
$$\Tr(S)=\sum\{X_p; \exists q \in P :X_q \cap S \text{ is stationary and }
 q<_P p \} $$
in the Boolean algebra $\Cal P (\Reg(\k))/\NS.$
Moreover for any $S\subseteq X_p$ stationary $o(S)=o(X_p)=o_P(p)$
where $o_P(p)$ is the rank of $p$ in $P$ and $o(\k)=o(P)=
\sup \{ o_P(p)+1; p\in P \} .$

In case  $P$ is linerly ordered 
and $S<\Reg(\k)$ for any stationary $S\subseteq\Sing(\k)$
we say that {\it the Axiom of Full Reflection
holds at} $\k$ [JS93], the sets $X_p$ are then the canonical stationary
sets (see [J84]).
[JS93] proves realization of well-orderings of length $\leq \k^+$, 
[JW93] even of length $>\k^+ .$

Note that only posets of cardinality $\leq\k^+$ can be realized as 
we have only $\k^+$ subsets of $\k$ (assuming GCH throughout the paper).

A sequence $S=\langle S_\lambda ; \lambda \leq \k \rangle$ is called
{\it a closed system of measures} (see [Mi83]) if
$$ \forall \lambda \leq \kappa \forall W\in S_\lambda : (j_W S)(\lambda)
   \subseteq S_\lambda $$
where $j_W$ is the canonical embedding $j_W:V\rightarrow V^\lambda /W .$
 For $U,W \in S_\lambda \; (\lambda \leq \k)$
define $U \tS W$ if $U\in (j_W S)(\lambda) .$ It is easy to prove
that $\tS$ is transitive. The standard Mitchell ordering extends $\tS$,
consequently $\tS$ must be well-founded.
Let $o_S (U)$ denote the rank of $U$ in this ordering.
 If $\tS$ is linear (on all 
$S_\lambda$) then the system is called {\it a coherent sequence of measures}
(see [Mi80],[JW93]). We say that the measures in $S_\k$ are {\it separable}
if there are sets $\langle X_U ; U \in S_\k \rangle$ so that 
$$ \forall U,W \in S_\k : \; X_U \in W \text{ iff } U=W .$$
If $S=\langle S_\lambda ; \lambda \leq \k \rangle$ is a 
closed system of measures
then we say that a $U \in S_\k$ is {\it a repeat point} if 
$$ \forall X\in U \exists W \tS U :\; X\in W $$
i.e. $U$ is not separable from its predecessors. It has been proved
in [JW93] that if
$S$ is a coherent sequence and 
 there are no repeat points in $S_\k$ then the
measures in $S_\k$ are separable.

In section 2 the construction of [JW93] is generalized to show the following

\proclaim{Theorem 1} If $S=\langle S_\lambda ; \lambda \leq \k \rangle$ is
a closed system of separable measures then there is a generic extension
$V[G]$ preserving cardinalities, cofinalities, and GCH where the reflection 
ordering of stationary subsets of $\Reg(\k)$ realizes the poset $(S_\k,
\tS)$ (as computed in $V$).
\endproclaim

Section 3 analyzes the question what well-founded posets are representable
by $(S_\k,\tS)$ and when the measures in $S_\k$ are separable.
It turns that closed systems of measures can be easily constructed
using a Laver's function on $\k$ that exists in $L[\Uarr]$ if
$\oU(\k)=\k^{++} .$

\proclaim {Theorem 2} Assume that $V=L[\Uarr],$ $\oU(\k)=\k^{++} ,$ and that
$P$ is a well-founded poset of cardinality $\leq\k^+ .$ Then there is a 
closed system of measures $\langle S_\lambda; \lambda\leq\k \rangle$
such that $P\cong (S_k,\tS )$ and the measures in $S_\k$ are separable.
\endproclaim

\proclaim {Corollary} If $V=L[\Uarr],$ $\oU(\k)=\k^{++}, $ then any
well-founded poset of cardinality $\leq\k^+$ is realized by the reflection 
ordering of stationary subsets of $\Reg(\k)$
in a generic extension of $V$ preserving cofinalities, 
cardinalities, and GCH.
\endproclaim

\vskip8pt
\subhead
2. The forcing construction
\endsubhead
\vskip8pt

The construction is analogous to the construction of [JW93]. 
We will spell out its definition but  will not repeat the
proofs that are almost literally same as the proofs in [JW93].

Let $S=\langle S_\lambda ; \lambda \leq \k \rangle$ be a closed
system of measures in the ground model $V$ satisfying GCH.

As usual, if $P$ is a forcing notion then $V(P)$ denotes either
the Boolean valued model or a generic extension by a $P$-generic 
filter over  $V$.

$P_{\k+1}$ will be an Easton support iteration of $\langle Q_{\lambda};\; 
\lambda \leq\k \rangle$, $Q_{\lambda}$
will be nontrivial only for $\lambda$ Mahlo.
$Q_{\lambda}$ (for $\lambda$ Mahlo) is defined in $V(P_{\lambda})$,
where $P_{\lambda}$ denotes the iteration below $\lambda$, 
as an iteration of length
$\lambda^{+}$ with $<\lambda$-support of forcing notions shooting clubs
through certain sets $X\subseteq \lambda$ (we will denote this standard forcing
notion $CU(X)$), always with the property that $X\supseteq \Sing(\lambda)$.
This condition will guarantee $Q_{\lambda}$ to be essentially $<\lambda$-closed
(i.e. for any $\gamma < \lambda$ there is a dense $\gamma$-closed subset of 
$Q_{\lambda}$).
$Q_{\lambda}$ will also satisfy the $\lambda^{+}$-chain condition. Consequently
$P_{\lambda}$ will satisfy $\lambda$-c.c. and will have size $\lambda$. 
Cardinalities, cofinalities, and GCH will be preserved, stationary subsets of 
$\lambda$
can be made nonstationary only by the forcing at $\lambda$, not below $\lambda$,
and not after the stage $\lambda$ - after stage $\lambda$ no subsets of 
$\lambda$ are added.

We use the $\lambda^{+}$-chain condition  of $Q_{\lambda}$ to get a canonical
enumeration of length $\lambda^+$ of all  $Q_{\lambda}$-names 
for subsets of $\lambda$ so that the $\beta$th name appears in $V(P_{\lambda} 
\ast Q_{\lambda} | \beta)$. Moreover for $U \in S_\lambda$ we will 
define certain filters $F_U$ in $V(P_\lambda \ast Q_\lambda 
| \beta)$.
Their definition will not be absolute, however the filter $F_U$ will extend the
measure $U$ and will 
increase coherently during the iteration.         

\proclaim{Definition} An iteration $Q$ of $\langle CU(B_\a); \a<\a_0 \rangle$
with $<\lambda$-support and length $\a<\lambda^+$ is called {\it an 
$\tilde S$-iteration} (for $\tilde S \subseteq S_\lambda$) if for all
$\a<\a_0$
$$V(P_\lambda \ast Q|\a) \models \text{ "$B_\a \in F_U$ for any $U \in
   \tilde S$ and $\Sing(\lambda) \subseteq B_\a$."}  $$
An $(S_\lambda|U \cup \{U \})$-iteration, where $S_\lambda |U = \{
W \in S_\lambda ; W \tS U \} ,$ is also called 
{\it an iteration of order $U .$}
 \endproclaim

(Note that any $S_\lambda$-iteration is also an $S_\lambda |U$-iteration.)

$Q_{\lambda}$ is then defined as an iteration of
$\langle CU(B_\alpha);\; \alpha < \lambda^{+} \rangle$
with $<\lambda$-support and length $\lambda^{+}$ so that every
$Q_\lambda | \alpha$ is an $S_\lambda$-iteration 
 and all potential names  $\dot X \subseteq \lambda$
are used cofinally many times in the iteration as some $B_\alpha$.

Observe that $Q_\lambda$ can be represented in
$V(P_\lambda)$ as a set of sequences of closed bounded subsets of
$\lambda$ in $V(P_\lambda)$ rather than in $V(P_\lambda \ast Q_\lambda 
| \alpha)$.
Moreover if $\dot q$ is a $P_\lambda$-name such that $1 \force_{P_\lambda}
\dot q \in Q_\lambda$ then using the $\lambda$-chain condition of $P_\lambda$
there is a set $A \subseteq \lambda^{+}$ (in $V$) of cardinality $<\lambda$
and $\gamma_0 < \lambda$ so that 
$$1\force_{P_\lambda} \;\; \text{supp}\;\dot q
\subseteq A  \text{ and } \forall \alpha \in A : \; \dot q (\alpha)
\subseteq \gamma_0 \;\;.$$
Consequently, $Q_\lambda$ can be represented as a set of functions
$g:A \times \gamma_0 \rightarrow [P_\lambda]^{<\lambda}$ where 
$A \subseteq \lambda^{+}$, $|A|<\lambda$ and $\gamma_0 < \lambda$.
In this sense $Q_\lambda$ has cardinality $\lambda^{+}$ and any
$Q_\lambda|\alpha$ has cardinality at most $\lambda$.

\vskip8pt
\subhead
Definition of filters $F_U$
\endsubhead
\vskip8pt

The filters $F_U$ ($U\in S_\lambda$) are defined 
in $V(P_\lambda \ast Q)$, where $Q$ is any iteration of 
order $U$, by induction so that the following is satisfied:

\proclaim{Proposition 2.1} Let $Q$, $Q'= Q\ast R$ be two iterations
of order $U$ then
$$F_{U}^{V(P_\lambda \ast Q)} \; = \; F_{U}^{V(P_\lambda
\ast Q')} \; \cap \; V(P_\lambda \ast Q) \; \; .$$
Moreover  $\; F_{U}^{V(P_\lambda)}\cap V = U$.
\endproclaim

\proclaim{Proposition 2.2} Let $j=j_{U}$ be the canonical embedding
from $V$ into $V^{\lambda}/U = M$ and
$Q$  an iteration of order $U$. Then $j$ can be lifted to an elementary 
embedding from a generic extension $V(P_\lambda \ast Q)$ of $V$ to
a generic extension $M(jP_\lambda \ast jQ)$ of $M$.
\endproclaim

\proclaim{Lemma 2.3} Let $ N=V^{\lambda}/W$ for some
$W\vartriangleright_S U$ and $Q$ be an iteration of order $U$. Then
$$F_{U}^{V(P_\lambda \ast Q)} \; = \; F_{U}^{N(P_\lambda \ast Q)}
\; \; . $$
Note that it also means that the definition of $F_{U}$
relativized to $N(P_\lambda \ast Q )$ makes sense.
\endproclaim

\proclaim{Lemma 2.4} Let $j=j_{U}:V \rightarrow M$. Then any 
$S_\lambda|U$-iteration $Q$ 
 is an subiteration of $(jP_\lambda)^{\lambda}$,
where $(jP_\lambda)^{\lambda}$ is the factor of
$jP_\lambda = P_\lambda \ast (jP_\lambda)^{\lambda} \ast 
(jP_\lambda)^{>\lambda}.$ Consequently
for any $G^{*}$ $jP_\lambda$-generic$/V$ and any $q \in Q$ there is an
$H \in M[G^{*}]$ $Q$-generic$/V[G]$ containing $q$
given by an embedding of $Q$ as a subiteration of $(jP_\lambda)^{\lambda}$, 
where $G=G^{*} \restriction
P_\lambda$.
\endproclaim

\proclaim{Definition}
Let $j,\; Q,\; G^{*},\; G$ be as in the lemma. Then
$\Gen_j (Q,G^{*})$ denotes the set of all filters $H\in M[G^{*}]$ 
$Q$-generic$/V[G]$ given by an embedding of $Q$ as 
a subiteration of $(jP_\lambda)^{\lambda}.$
\endproclaim

\proclaim{Lemma 2.5} Let $j$ be as above, $Q$ an iteration of order
$U$, $G^{*}$ $jP_\lambda$-generic$/V$, $H \in \Gen_j (Q,G^{*})$.
For every $\beta < l(Q)$ let $C_\beta \subset \lambda$ be the club
$\cup\{r(\beta); r\in H\}$, and let
$[H]^{j}$ denote the $j(l(Q))$-sequence given by
$$[H]^{j}(\gamma) = \cases C_\beta\cup \{ \lambda \} , & 
\text{ if } \gamma=j(\beta) \\
                           \emptyset & \text{ otherwise.} \endcases
$$
Then $[H]^{j} \in jQ_{/G^{*}}$.
\endproclaim

Propositions 2.1 and 2.2 are essential to analyze reflection of stationary
sets
in the generic extension.

Now let $U\in S_\lambda$ and suppose that the filters $F_{U'}$ 
have been defined for all $\lambda'< \lambda$ and $U' \in S_{\lambda'}$
and for all $U' \in S_\lambda |U$ so that
2.1-2.5 holds. Moreover
let $\alpha < \lambda^{+}$ and $F_U$ be defined for 
all iterations
of order $U$ and length $<\alpha$ so that 2.4 and 2.5 holds for $U$
and iterations of length $\leq \alpha$. Then we can define $F_U$
for iterations of order $U$ and length $\alpha$.

\proclaim{Definition} Let Q be an iteration of order $U$ and length
$\alpha$, $j=j_U :V \rightarrow M.$ For a $P_\lambda \ast Q$-name 
$\dot X$ of a subset of $\lambda$
and $p\ast q \in P_\lambda \ast Q$ define
$$ p\ast q \force _{P_\lambda \ast Q} \; \dot X \in F_\delta $$
if the following holds in V:

$$j(p) \force_{jP_\lambda}\text{ ``For any $H \in \Gen_j (Q,G^{*})$ containing 
$q$ :} $$
$$[H]^{j} \force_{jQ} \; \check \lambda \in j\dot X \;\;\; \text{.''}$$
\endproclaim

The definition says that $p\ast q \force \dot X \in F_\delta$ if 
$\lambda \in j^{*} X$ whenever $j^{*}:V[G\ast H] \rightarrow M[G^{*} 
\ast H^{*}]$
is a lifting of $j$ of certain kind and $p\ast q\in G\ast H$.

For the proofs of 2.1--2.5 see [JW93], modifications are
left to the reader.

\vskip8pt
\subhead
Reflection of stationary sets in the generic extension
\endsubhead
\vskip8pt

Let us analyze the behaviour of the reflection ordering of stationary subsets
of $\lambda$ in the generic extension $V(P_{\lambda+1})$.
That is the same as in $V(P_{\k+1})$ because no subsets of $\lambda$
are added after the stage $\lambda .$
Note that any subset of $\lambda$ in $V(P_{\lambda+1})$ already
appears in $V(P_\lambda \ast Q_\lambda|\a)$ for some $\a<\lambda^+ .$

\proclaim{Lemma 2.6}
$V(P_\lambda \ast Q_\lambda |\alpha) \models \Club(\lambda) \subseteq F_U$
for any $U \in S_\lambda. $
\endproclaim

\demo{Proof}
See 3.2 in [JW93].
\enddemo

\proclaim{Lemma 2.7} Let $S\subseteq \Reg(\lambda)$ be in $V(P_{\lambda+1})$
then $S$ is stationary iff there is a $U \in S_\lambda$ such that for some
$\a < \lambda^+$ $S \in V(P_\lambda \ast Q_\lambda|\a)$ and 
$V(P_\lambda \ast Q_\lambda|\a) \models$ ``$S$ is $F_U$-positive''
(i.e. iff this is true for all $\a$ such that $S \in V(P_\lambda \ast
Q_\lambda|\a)$).
\endproclaim

\demo{Proof}
Let $S \in V(P_\lambda \ast Q_\lambda|\a)$ and
$V(P_\lambda \ast Q_\lambda|\a) \models$ ``$S$ is $F_U$-positive''
and suppose by contradiction that $V(P_{\lambda+1}) \models$ ``$S$ is
nonstationary''.
Then there is an $\a' > \a$ such that
$$V(P_\lambda \ast Q_\lambda|\a') \models \text{``$S$ is nonstationary''}$$
however proposition 2.1 says that
$$V(P_\lambda \ast Q_\lambda|\a') \models \text{``$S$ is $F_U$-positive''}$$
and that is a contradiction with lemma 2.6.

On the other hand suppose that $S \in V(P_\lambda \ast Q_\lambda|\a)$
and 
$$V(P_\lambda \ast Q_\lambda|\a) \models \text{``$S$ is $F_U$-thin for any
   $U \in S_\lambda$ .''}$$
Then by the definition of $Q_\lambda$ a club is shot through $\k \setminus S$,
consequently $S$ is nonstationary in $V(P_{\lambda +1}).$
\qed
\enddemo

\proclaim{Lemma 2.8} Assume that $S$ is in $V(P_\lambda \ast Q_\lambda |\a),$ 
$S\subseteq \Reg(\lambda)$ and $W\in S_\lambda .$
Then $\Tr(S)$ is either $F_W$-thin or it belongs to $F_W .$
The former is true iff $S$ is $F_U$-thin for all $U \tS W .$
Moreover if $S\subseteq \Sing(\lambda)$ is stationary
then $\Tr(S) \in F_U$ for all $U\in S_\lambda .$
\endproclaim

\demo{Proof}
Suppose that there is $U \tS W$ such that
$p\ast q \force _{P_\lambda\ast Q_\lambda |\alpha} \; \text{``$\dot S$ 
is $F_U$-positive''}$ but
$p\ast q \notforce _{P_\lambda\ast Q_\lambda |\alpha} \;$``$\Tr(\dot S) \in 
F_{W}$'' for some $W \vartriangleright_S U.$
Denote $j=j_W :V \rightarrow M=V^{\lambda}/W .$
Then there is a filter $G^{*}$ $jP_\lambda$-generic$/V,$ $G^{*} \ni jp,$
and $H \in  \Gen_j (Q_\lambda|\alpha,G^{*}),$ $H \ni q,$
and a filter $H^{*}$ $j(Q_\lambda |\alpha)$-generic$/V[G^{*}],$
$H^{*}$ contains $[H]^{j} ,$
so that $\lambda \notin j\Tr(\dot S)_{/G^{*} \ast H^{*}}.$
The embedding $j$ is lifted to
$$ j^{**}: V[G\ast H] \rightarrow M[G^{*} \ast H^{*}] \; . $$
$S=\dot S_{/G\ast H}$ is $F_U$-positive in $V[G\ast H]$ and
$$j\Tr(\dot S)_{/G^{*} \ast H^{*}} \; = \; j^{**}(\Tr(S)) \; = \; 
\Tr^{M[G^{*}\ast H^{*}]}(j^{**}S)\; .$$
Thus $\lambda \notin \Tr^{M[G^{*}\ast H^{*}]}(j^{**}S)$ which means that
$$ M[G^{*} \ast H^{*}] \models \text{`` $S$ is not stationary in $\lambda$''}$$
because $S=j^{**}S \cap \lambda.$ Consequently 
$$V[G^{*}\ast H^{*}] \models  \text{``$S$ is not stationary in $\lambda$.''}$$
Observe that $j(Q_\lambda |\alpha)$ has a dense subset $\lambda$-closed
in $M[G^{*}]$ and thus also in $V[G^{*}].$ Moreover
$jP_\lambda = P_\lambda \ast (jP_\lambda)^{\lambda} \ast R$ where $R$
is essentially $\lambda$-closed in $V[G^{*} | \lambda +1].$ It implies that
already $V[G^{*} | \lambda +1] \models  \text{``$S$ is not stationary 
in $\lambda$.''}$
Let us now consider the isomorphism $(jP_\lambda)^{\lambda} \simeq
(Q_\lambda |\alpha) \ast \tilde Q$ from the proof of 2.4 giving
a filter $H=G^{*} \restriction (Q_\lambda |\alpha),$
let $\tilde H = G^{*} \restriction \tilde Q.$ Since every subset of $\lambda$
in $V[G\ast H \ast \tilde H]$ is already in some $V[G\ast H \ast 
\tilde H|\beta]$ there is a $\beta < \lambda^{+}$ so that
$$V[G\ast H \ast \tilde H|\beta] \models \text{``$S$ is not stationary 
in $\lambda$.''}$$
But since $(Q_\lambda |\alpha) \ast (\tilde Q |\beta)$ is an 
$S_\lambda|W$-iteration and hence an iteration of order $U$
it follows from proposition 2.1 that
$$V[G\ast H \ast \tilde H|\beta] \models \text{``$S$ is $F_U$-positive''}$$
which contradicts 
lemma 2.6.

The proof that for any stationary $S\subseteq \Sing(\lambda)$
and $U \in S_\lambda,$ $\Tr(S) \in F_U$ is analogous using the following fact instead of proposition 2.1.
 
\proclaim{Claim}
Stationary subsets of $\Sing(\lambda)$ are preserved by $\emptyset$-iterations.
\endproclaim
\demo{Proof} 
See 7.38 in [J86] or 3.4 in [JW93].
\enddemo

Now let $S \in V(P_\lambda \ast Q_\lambda |\a)$, $W \in S_\lambda$ and
$S \subseteq \Reg(\lambda)$ be $F_U$-thin for any $U \tS W .$ We want to prove 
that $V(P_\lambda \ast Q_\lambda |\a) \models \lambda \setminus \Tr(S) \in F_W .$
Let $j=j_W:V \rightarrow M,$  then $(jP_\lambda)^\lambda$ is an iteration
of length $\lambda^+$ such that $(jP_\lambda)^\lambda|\beta$ is always
an $S_\lambda|W$-iteration and every potential name is used cofinally
many times. For $\beta$ large enough $Q_\lambda|\a$ is an subiteration
of $(jP_\lambda)^\lambda|\beta$ and
$$V(P_\lambda \ast (jP_\lambda)^\lambda |\beta) \models 
\text{``$\lambda\setminus S
\in F_U$ for any $U \in S_\lambda |W$ ''}$$
because $F_U^{ V(P_\lambda \ast (jP_\lambda)^\lambda |\beta)} 
\supseteq F_U^{V(P_\lambda \ast Q_\lambda|\a)} \ni \lambda\setminus S .$
Consequently a club is shot through $\lambda\setminus S$ in the
iteration $(jP_\lambda)^\lambda.$ It implies that
$$V[G^*] \models \text{``$S\subseteq \lambda$ is nonstationary''}$$
where $G^*$ is any $jP_\lambda$-generic/$V$ and thus also
for any $H \in \Gen_j(Q_\lambda, G^*)$ and any $H^* \ni [H]^j$
$j(Q_\lambda |\a)$-generic/$V[G^*]$
$$V[G^*\ast H^*] \models \lambda \in j^*(\lambda\setminus \Tr(S)) .$$
That proves 
 $V(P_\lambda \ast Q_\lambda |\a) \models \lambda \setminus \Tr(S) \in F_W .$
\qed
\enddemo

The filters $F_U$ are not defined in $V(P_\lambda \ast Q_\lambda)$,
however we can define $\tilde F_U = \bigcup_{\a<\lambda^+} F_U^{V(P_\lambda
\ast Q_\lambda|\a)} .$
It follows from the lemmas that a set $S\subseteq \Reg(\lambda)$
in $V(P_{\lambda +1})$ is stationary iff it is $\tilde F_U$-stationary
for some $U.$ Moreover for any $U$ either $\Tr(S) \in \tilde F_U$ or
$\Tr(S)$ is $\tilde F_U$-thin . The former is true iff there is a $W \tS U$
such that $S$ is $\tilde F_W$-positive. Besides for 
$S\subseteq \Sing(\lambda)$ stationary $\Reg(\lambda) \subseteq \Tr(S)$
(mod $\NS$).
If the measures $U \in S_\lambda$ are separated by  sets 
$X_U \subseteq \Reg(\lambda)$  it follows that
$\langle X_U; U\in S_\lambda \rangle $ forms a maximal antichain
of stationary subsets of $\Reg(\lambda)$ in $V(P_{\lambda+1})$
and the reflection ordering of stationary subsets of $\Reg(\lambda)$
realizes the poset $(S_\lambda, \tS) .$
Actually $\tilde F_U=C[X_U]=
\{Y\subseteq \k ; \exists C\text{ a club}:C\cap X_U \subseteq Y \} .$
That proves theorem ~1.

\vskip8pt
\subhead
3. Closed systems of measures
\endsubhead
\vskip8pt

In this section we construct closed systems of separable measures
isomorphic to a given well-founded poset.

\proclaim{Lemma 3.1}
Let $S=\langle S_\lambda ; \lambda \leq \k \rangle $ be a system of measures
such that 
$$\forall U,W \in S_\k :\; (j_U S)(\k) \subseteq S_\k 
  \text{ and }
  (W\in (j_U S)(\k) \Rightarrow (j_W S)(\k)\subseteq (j_U S)(\k)). $$
Then there is a closed subsystem of measures 
$\tilde S=\langle \tilde S_\lambda ; \lambda \leq \k \rangle $
such that $S_\k = \tilde S_\k$ and
$\forall U,V \in S_\k: U \tS V \leftrightarrow U\vartriangleleft_{\tilde S} V .$
\endproclaim

\demo{Proof}
Define $\tilde S_\lambda$ by induction on $\lambda \leq \k .$
Suppose $\tilde S_{\lambda'}$ has been defined for $\lambda'<\lambda .$
If $\forall U\in S_\lambda : (j_U (\tilde S\restriction \lambda))(\lambda)
 \subseteq S_\lambda$
then put $\tilde S_\lambda = S_\lambda$ otherwise 
$\tilde S_\lambda = \emptyset .$
That defines $\tilde S_\lambda \subseteq S_\lambda$ for all $\lambda \leq \k ,$
obviously $\tilde S_\k = S_\k$ by the assumption of the lemma.
Finally let us prove that for any $U\in S_\k$
$$ \{ \lambda < \k;\; \tilde S_\lambda = S_\lambda \} \in U $$
and consequently $(j_U S)(\k) = (j_U \tilde S)(\k) .$
It is enough to prove that the following set is in $U$
$$ A = \{ \lambda < \k; \forall W \in S_\lambda: (j_W (S\restriction
 \lambda))(\lambda)
 \subseteq
   S_\lambda \} .$$
That is true iff $\k \in j_U A,$ iff
$$ \forall W \in (j_U S)(\k) : (j_W ^{M_U}(S\restriction \k))(\k) 
 \subseteq (j_U S)(\k) $$
where $j_U:V \rightarrow V^\k /U = M_U$ and
$j_W ^{M_U} : M_U \rightarrow M_U^\k /W .$
For any $W \in (j_U S)(\k) \subseteq S_\k$ the embedding 
$j_W ^{M_U} = j_W \restriction  M_U , $ thus we have
$(j_W ^{M_U}(S\restriction \k)(\k) = (j_W S)(\k) \subseteq (j_U S)(\k) .$
\qed
\enddemo

As a corollary we can prove

\proclaim{Proposition 3.2}
Suppose $P$ is a  well-founded  poset of cardinality $\leq \k$
and  there are measures $U_p$ ($p\in P$) on $\k$ so that
$U_p \vartriangleleft U_q$ whenever $p <_P q .$
Then there is a closed system of separable measures
$S=\langle S_\lambda ; \lambda \leq \k \rangle $ such that
$S_\k = \{ U_p ; p\in P \}$ and
$$p<_P q \leftrightarrow U_p \tS U_q .$$
\endproclaim

\demo{Proof}
Since $|P| \leq \k$ we can find disjoint sets $X_p \in U_p$ ($p\in P$)
separating the measures. For $p\in P$ the set $\{ U_q ; q<_P p \}$
is in the ultraproduct $V^\k / U_p$ because all $U_q$ ($q<_P p$)
are and  the number of them is at most $\k .$
Let $\langle S_\lambda ; \lambda \in X_p \rangle$ be a sequence 
of sets of measures over $\lambda$ such that
$$ [S_\lambda]^{\lambda \in X_p}_{U_p} = \{ U_q; q<_P p \} . $$
Glue together all those sequences into 
$\langle S_\lambda ;\lambda \leq \k \rangle$ (put $S_\lambda = \emptyset$
if $\lambda \notin \bigcup _{p\in P} X_p$ ). Then for any $p\in P$
$$(j_{U_p}S)(\k) = \{ U_q; q <_P p \} .$$
Finally use lemma 3.1 
(the assumption is satisfied as $<_P$ is transitive)
to get the desired closed system of measures.
\qed
\enddemo

Closed systems of measures can be easily constructed using a Laver's
function.

\proclaim{Definition} {\rm (cf. [La78])} We say that
a function $f:\k\ra V_\k$ is 
{\it a Laver's function on} $\k$
 if
$$\forall x\in V_{\k+2}\exists U \text{ a measure on }\k:(j_Uf)(\k)=x .$$
\endproclaim

It means that the Laver's function serves as a universal function
in $\;^\k V$ for all
$x\in V_{\k+2}.$ 
Generalizing [La78] we can prove that if $\k$ is $\Cal P_2 \k$-strong
then there is a Laver's function on $\k .$
We will show that the existence of a Laver's function
on $\k$ is 
actually equiconsistent with the Mitchell order of $\k$ being $\k^{++} .$

It will follow that if there is
a Laver's function on $\k$ then 
there is a coherent sequence of measures $\Uarr$ such that 
$o^{\Cal U}(\k)=\k^{++},$ and  moreover that {\it the measures
on $\k$ cover $\Cal P(\k^+)$} in the following sense:
$$\forall A\in \Cal P(\k^+) \exists \a<\k^{++} : A\in \Ult(V,U_\a^\k) .$$
Let us  prove that those two conditions are sufficient for the existence
of a Laver's function.

\proclaim{Proposition 3.3} Let $\Uarr$ be a coherent sequence of measures
such that $o^{\Cal U}(\k)=\k^{++}$ and the measures 
$\{U_\a^\k; \a<\k^{++}\}$ cover $\Cal P(\k^+) .$ Then there is 
a Laver's function on $\k.$
\endproclaim

\demo{Proof} Firstly observe that $\{U_\a^\k; \a<\k^{++}\}$ also
covers $V_{\k+2}:$ let $\pi\in \Ult(V,U_0^\k)$ be a bijection between
$\k^+$ and $V_{\k+1} .$ If $x\in V_{\k+2} ,$ put $A=\pi^{-1}[x] ,$
and find $\a<\k^{++}$ such that $A\in \Ult(V,U_\a^\k) .$
Since $\pi\in \Ult(V,U_\a^\k)$, $x$ must be in $\Ult(V,U_\a^\k).$
It can be  assumed without loss of generality that
$\oU(\l)<\l^{++}$ for all $\l<\k .$

\proclaim{Claim} There are well-orderings $R^\l$ of $V_{\l+2}$
of order type $\l^{++}$
($\l\leq \k$) so that for any $\a<o^{\Cal U}(\k)$ the well-ordering
$j_{U_\a^\k}(\langle R^\l;\l<\k\rangle)(\k)$
is an initial segment of $R^\k .$
\endproclaim

\demo{Proof} Start with arbitrary well-orderings $<_\l$ of $V_{\l+2}$
of order type $\l^{++}$
($\l\leq\k$). Assume $\langle R^\l;\l<\mu \rangle$ have been defined
so that for all $\l<\mu$
$$\forall\a<\oU(\l) : R^\l_\a=j_{U_\a^\l}(\langle R^{\l'};\l'<\l\rangle)(\l)
   \text{ is an initial segment of }R^\l . \tag3.1 $$
Let $\beta<\oU(\mu),$ then by elemetarity
 $R^\mu_\beta=
   j_{U^\mu_\beta}(\langle R_\l;\l<\mu \rangle)(\mu)$
is a well-ordering of $V_{\l+2}\cap\Ult(V,U_\beta^\mu)$
of order type $\mu^{++\Ult(V,U_\beta^\mu)}$ and
$$\Ult(V,U^\mu_\beta)\models \lqot \forall \a<\beta :
  R^\mu_\a=j_{U^\mu_\a}(\langle R_\l;\l<\mu\rangle)(\mu)
  \text{ is an initial segment of }R^\mu_\beta . \rqot $$
Consequently $R^\mu$ can be defined as $\cup\{R^\mu_\beta; \beta<\oU(\mu)\}$
and end-extended by $<_\mu$ on the remaining elements of $V_{\mu+2} .$
Then (3.1) is also satisfied for $\l=\mu .$ If $\mu=\k$ then there are no
remaining sets in $V_{\k+2},$ and $R^\k=\cup\{R^\k_\beta; \beta<\oU(\k)\} .$
\qed $\;$Claim
\enddemo
Now define $f$ by induction as follows: assume $f(\l')$ has been defined
for $\l'<\l ,$ and let $f(\l)$ be the $R^\l$-least $x\in V_{\l+2}$
such that
$$\neg \exists\a<\oU(\l): x=(j_{U^\l_\a}(f\res\l))(\l) .$$
Leave $f(\l)$ undefined if there is no such $x.$ Assume there  is
$x\in V_{\k+2}$ such that
$$\neg \exists\a<\oU(\k) : x=(j_{U^\k_\a} f)(\k) ,$$
and let $x$ be the $R^\k$-least with this property.
Let $\beta<\oU(\k)$ be such that $x \in \Ult(V,U^\k_\beta)$ and
$$\forall y <_{R^\k} x \exists \a<\beta : y=(j_{U^\k_\a} f)(\k) ,$$
this is possible as $\oU(\k)=\k^{++}$ and there are at most $\k^+$ many
$y<_{R^\k} x .$ Then
$$ \Ult(V,U^\k_\beta)\models \lqot x\text{ is the $R_\beta^\k$-least such that }
\neg\exists\a<o^{j\Cal U}(\k)=\beta : 
   x=(j_{U^\k_\a} f)(\k). \rqot $$
Consequently $(j_{U^\k_\beta} f)(\k)=x$ by elementarity
- a contradiction. \qed $\;$3.3
\enddemo

Let us prove that the assumptions of 3.3 are much weaker than 
$\Cal P_2\k$-strongness. W.Mitchell proved in [Mi83] that there
is an inner model $L[\Uarr]$ satisfying GCH such that $\Uarr$ is a 
coherent sequence of measures in $L[\Uarr],$ all measures in $L[\Uarr]$
are in $\Uarr,$ and $\oU(\k)=\min(o(\k),\k^{++L[\Uarr]})$ for all ordinals
$\k .$

\proclaim{Proposition 3.4} If $\oU(\k)=\k^{++}$ in $L[\Uarr]$ then
the measures on $\k$ cover $\Cal P(\k^+) .$
\endproclaim

\demo{Proof} Let us firstly suppose that $\k$ is the maximal measurable
cardinal in $L[\Uarr] $ (it can be easily achieved by cutting the universe
at the first measurable above $\k$ and then applying the Mitchel's
construction of $L[\Uarr]$). Notice that the sequence of measures
must be actually represented as 
$\Uarr = \{ ( \l,\a,X ) ; X\in U^\l_\a \}$
in order $L[\Uarr]$ makes sense. Let $A\subseteq \k^+,$ 
$A\in L_\gamma[\Uarr],$ $\gamma >\k^+ .$ Find a model 
$M\subseteq L_\gamma[\Uarr]$ of cardinality $\k^+$ such that
$\k^+ +1 \subseteq M,$ $ A\in M,$ $ \Cal P(\k)\subseteq M,$
$$\forall \a\in M : \a<\k^{++} \Rightarrow \a\subseteq M ,\text{ and} $$ $$
\langle M,\in, \Uarr \cap M \rangle \prec \langle L_\gamma[\Uarr],\in,\Uarr
                                          \rangle .$$
Let $\pi$ be the transitive collapse of $M .$
Then it is easy to see that 
$\overrightarrow{\Cal V}=\pi(\Uarr\cap M)=\Uarr\cap M=\Uarr\res(\k,\Theta)$
where $\Theta=M\cap\k^{++} .$ 
By 32.7 of [Ka93] $\pi[M]=L_\delta[\Uarr\res(\k,\Theta)]$ for a 
$\delta<\k^{++} .$ Consequently
$$A=\pi(A)\in L_\delta[\Uarr\res(\k,\Theta)]\subseteq \Ult(V,U^\k_\Theta).$$

The idea of the proof in general situation is due to W. Mitchell
(personal communication).
Start with $A\subseteq\k^+,$ let $M,$ $L_\gamma[\Uarr],$ $\pi,$
$\overrightarrow{\Cal V}=\pi(\Uarr\cap M),$
$L_\delta[\overrightarrow{\Cal V}]=\pi[M]$ be as above.
Then we can prove only that 
$\overrightarrow{\Cal V}\res(\k+1) =\Uarr\res (\k,\Theta),$
the measures above $\k$ in $\overrightarrow{\Cal V}$ does not have 
to be same as in $\Uarr .$
Observe that $L_\delta[\overrightarrow{\Cal V}]$ is iterable since it is 
embedded into iterable $L_\gamma[\Uarr] .$
By the Mitchell's comparison lemma [Mi83] there are iterations
defined in $L[\Uarr]$
$$i:L_\delta[\overrightarrow{\Cal V}]\ra L_\vartheta[\overrightarrow{\Cal W}],$$
$$i':L[\Uarr]\ra L[\overrightarrow{\Cal W}']$$
such that either $\overrightarrow{\Cal W}$ is an initial segment of
$\overrightarrow{\Cal W}'$ or $\overrightarrow{\Cal W}'$ is an initial
segment of $\overrightarrow{\Cal W}.$
$\overrightarrow{\Cal W}$ is an initial segment of
$\overrightarrow{\Cal W}'$ if
$\overrightarrow{\Cal W}\cap L_\vartheta[\overrightarrow{\Cal W}]
\cap L[\overrightarrow{\Cal W}'] =
\overrightarrow{\Cal W}'\res \beta \cap L_\vartheta[\overrightarrow{\Cal W}]
\cap L[\overrightarrow{\Cal W}']$
for an ordinal $\beta .$
Assume  this is true, it is then easy to see that
$L_\vartheta[\overrightarrow{\Cal W}]\subseteq L[\overrightarrow{\Cal W}'] .$
It follows from the proof of the  comparison lemma that
$\crit(i)>\k,$ and the iteration $i'$ starts with the ultraproduct
by $U^\k_\Theta$ and then proceeds with measures above $\k .$
Thus $A\in L_\vartheta[\overrightarrow{\Cal W}]\subseteq
L[\overrightarrow{\Cal W}']$ implies
$A\in \Ult(L[\overrightarrow{\Cal W}], U^\k_\Theta ).$
In that case we are done.

Assume towards a contradiction that $\overrightarrow{\Cal W}'$
is a proper initial segment of $\overrightarrow{\Cal W},$ i.e.
$\overrightarrow{\Cal W}' \cap L_\vartheta[\overrightarrow{\Cal W}]
\cap L[\overrightarrow{\Cal W}'] =
\overrightarrow{\Cal W}\res \beta \cap L_\vartheta[\overrightarrow{\Cal W}]
\cap L[\overrightarrow{\Cal W}']$
and there is a measure in $\overrightarrow{\Cal W}$ above $\beta .$
It easily follows that $L_\vartheta[\overrightarrow{\Cal W}']=
L_\vartheta[\overrightarrow{\Cal W}\res\beta] .$
The top measure in $\overrightarrow{\Cal W}$
can be iterated making $\vartheta$ arbitrarily large,
hence $L[\overrightarrow{\Cal W}']=L[\overrightarrow{\Cal W}\res\beta]$
and $\overrightarrow{\Cal W}'=\overrightarrow{\Cal W}\res\beta 
\cap L[\overrightarrow{\Cal W}\res\beta] .$
Moreover the iterated ultraproduct of $L_\vartheta[\overrightarrow{\Cal W}]$
by the top measure produces a class of indiscernibles
(defined in $L[\Uarr]$) for $L[\overrightarrow{\Cal W}']$
containing all large enough cardinals.
That gives a definition in $L[\Uarr]$ of truth in 
$L[\overrightarrow{\Cal W}']$:
$L[\overrightarrow{\Cal W}']\models\sigma$ iff for all
regular cardinals $\l$ in $L[\Uarr]$ large enough
$V_\lambda \cap L[\overrightarrow{\Cal W}']\models\sigma .$
Consequently using the elementary embedding $i'$ we obtain 
a definition of truth in $L[\Uarr]$:
$L[\Uarr]\models\sigma$ iff for all
regular cardinals $\l$  large enough
$V_\lambda \cap L[\Uarr]\models\sigma .$
A contradiction. 
\qed
\enddemo

{\bf Remark} Notice that if $\Uarr$ is a coherent sequence of measures
such that $\oU(\k)<\k^{++}$ then the measures 
$\{U^\k_\a;\a<\oU(\k)\}$ cannot cover $\Cal P(\k^+):$ let $\gamma<\k^{++}$
be such that $\gamma > \k^{++\Ult(V,U_\a^\k)}$ for all $\a<\oU(\k) ,$
and let $A\subseteq \k^+$ code the well-ordering of order type $\gamma .$
Then $A\notin \Ult(V,U_\a^\k)$ for all $\a<\oU(\k).$ Hence
$\oU(\k)=\k^{++}$ is a necessary condition for the covering of 
$\Cal P(\k^+) .$ However not sufficient: Let $\Uarr$ be a coherent sequence
of measures in $V,$ $\oU(\k)=\k^{++}.$
Let $P$ be the Cohen forcing adding a subset of $\k^+,$
$G\subseteq\k^+$ $P$-generic over $V.$ Since $P$ is $\k$-closed
no new $\k$-sequences are added, and so $\Uarr$ is also a coherent sequence 
of measures in $V[G] .$ I claim that $G$ is not covered by any of the
measures $\{ U_\a^\k;\a<\k^{++}\}.$
Consider the ultraproduct embedding
$$j_\a^* :V[G]\ra \Ult(V[G],U_\a^\k)=M_\a[G^*]$$
which extends $j_\a:V\ra \Ult(V,U_\a^\k)=M_\a .$ $G^*=j_\a^*(G)$ is
$j_\a P$-generic over $M_\a .$ But $j_\a P$ is the Cohen forcing
adding a subset of $j_\a(\k^+)$ in the sense of $M_\a ,$ so it is
$j_\a\k$-closed in $M_\a ,$ in particular $\k^+$-closed, and no new
subsets of $\k^+$ are added by $G^* .$ As $G\notin M_\a$ we conclude that
$G\notin M_\a[G^*]=\Ult(V[G],U_\a^\k).$
It can be actually proved using Mitchell's methods that if $V=L[\Uarr]$ then
there are no new measures on $\k$ in $V[G] .$

Now let us use the Laver's function on $\k$ to construct closed
systems of measures isomorphic to a given well-founded poset. 

\proclaim {Proposition 3.5} If there is a Laver's function on $\k$
then there are two functions $F,G:\k\ra V_\k$ such that
$$\forall A,B\subseteq V_{\k+2}, |A|\leq\k^+, |B|\leq\k^+
  \exists U \text{ a measure on }\k:$$ $$
(j_U F)(\k)=A \text{ and } (j_U G)(\k)=B .$$
\endproclaim

\demo{Proof} All we need is an effective coding of pairs $(A,B),$
where $A,B\subseteq V_{\l+2},$ $|A|\leq\l^+,$ $|B|\leq\l^+,$
by elements of $V_{\l+2} .$
Firstly for $x,y\in V_{\l+2}$ define
$$x\oplus y = \{ \{0\}\times z;z\in x\} \cup \{ \{1\}\times z; z\in y\} .$$
It is easy to see that $x\oplus y \in V_{\lambda+2}$ if $\lambda$ is 
a limit ordinal. On the other hand for any $z\in V_{\lambda+2}$
 we can find unique $x,y\in V_{\lambda+2}$ such that $z=x\oplus y$
if there are any.
Given $A=\{x_i;i\in I\}\subseteq V_{\lambda+2}$ indexed over a set
$I\subseteq V_{\lambda+1}$ define
$$\bigoplus_{i\in I} x_i = \{i\times z; i\in I \text{ and }z\in x_i \}$$
which is in $V_{\lambda+2}$ if $\lambda$ is a limit ordinal. If
$z\in V_{\lambda+2}$ we can again find unique $I\subseteq V_{\lambda+1}$
and $A=\{x_i; i\in I\} \subseteq V_{\lambda+2}$ such that
$z=\bigoplus_{i\in I} x_i$ if there are any.
Finally code $(A,B),$ $A=\{x_i;i\in I\},$ $B=\{y_j;j\in J\}$ as
$$A\oplus B = (\bigoplus_{i\in I} x_i)\oplus (\bigoplus_{j\in J} y_j) .$$

Let $f:\k\ra V_\k$ be the Laver's function. For $\l<\k$ limit
put $F(\l)=A_\l,$ $G(\l)=B_\l$ if $f(\l)=A_\l\oplus B_\l$ for some
$A_\l,B_\l,$ otherwise $F(\l)=\emptyset=G(\l) .$ For a given pair
(A,B) find a measure $U$ on $\k$ such that $(j_Uf)(\k)=A\oplus B ,$
then by the construction $(j_U F)(\k)=A$ and $(j_U G)(\k)=B .$
\qed
\enddemo

\proclaim{Proposition 3.6}
Let $P$ be a  well-founded poset of cardinality $\leq \k^+ .$
Then there are two closed systems of measures 
$S=\langle S_\lambda;\lambda \leq\k\rangle ,$
$T=\langle T_\lambda;\lambda\leq\k\rangle$
such that $S_\k=T_\k$, $\vartriangleleft_T$ extends $\tS$,
$P\cong (S_\k,\tS)$ and $T$ forms a coherent sequence of measures
with $o_T(\k)\leq o(P)\cdot \k^+ .$ Moreover if $U_0$ is the first
measure in $T_\k$ then we can require that $j_{U_0}(\k)$ is greater
than a given ordinal $\vartheta < \k^{++} .$
\endproclaim

\demo{Proof} Let $F,G:\k\ra V_\k$ be the functions from 3.5.
Enumerate $P=\{p_\a;\a < \mu \}$ so that $p_a<_P p_\beta$ implies
$\a<\beta$ and each level of $P$ corresponds in this ordering
to a segment of  order type at most $\k^+ .$
Consequently $\mu\leq o(P)\cdot \k^+ .$
Firstly find $U_{p_0}$ such that $(j_{U_{p_0}}F)(\k)=\vartheta$
and $(j_{U_{p_0}}G)(\k)=\emptyset.$ Then $j_{U_{p_0}}(\k)>\vartheta$
because $(j_{U_{p_0}}F)(\k) \in V_{j_{U{p_0}}(\k)} .$
Using the two functions $F,G$ by induction on $\a>0$ find $U_{p_\a}$
so that
$$ (j_{U_{p_\a}}F)(\k)=\{U_{p_\gamma}; \gamma < \a \} \text{ and} $$
$$(j_{U_{p_\a}}G)(\k)=\{ U_{p_\gamma} ; p_\gamma <_P p_\a \} .$$
That gives $T_\k=S_\k=\{U_{p_\a} ; \a<\mu \} , $ put
$T_\lambda=F(\lambda)$ and $S_\lambda = G(\lambda)$ if those are
sets of measures over $\lambda .$
Use lemma 3.1 to get  closed systems of measures with
required properties.
\qed
\enddemo

The proof works for well-founded $\k^{++}$-like posets $P$ as well
(i.e. $|P|\leq\k^{++}$ and $\forall p\in P: |P\res p|\leq \k^+$).
 We will need the auxiliary 
coherent sequence of measures $T$ to make sure the measures in $S_\k$ are
separable.

The assumption $\oU(\k)=\k^{++}$ can be significantly
weakened to represent smaller well-founded posets.
The proofs of 3.3,3.6 can be modified to prove that if
 $\oU(\k)=\k^+$ then all well-founded $\k^+$-like posets $P$
(i.e. $|P|\leq\k^+,$ and $\forall p\in P: |P\res p|\leq \k$)
can be represented as $\langle S_\k,\tS \rangle$ for a closed system 
of measures $S.$

\vskip8pt
\subhead Separability and Repeat Points \endsubhead
\vskip8pt

To prove theorem 2 let us give some estimates on the order
of a least repeat point.

Fix a coherent sequence $\Uarr$ with the least repeat point on $\k$
of order $\Theta .$ 
We say that $f_\nu$ is {\it an $\a$-canonical function for $\nu$}
if for all $\delta \in [\a,\oU(\k)):$ $[f_\nu]_{U_\delta^\k}=\nu .$

\proclaim {Lemma 3.7} {\rm [JW93]} If $U_\a^\k$ is not a repeat
 point then $\a$
has an $\a$-canonical function. Moreover there is $X_\a\subseteq \k$
such that 
$$\forall\delta<\oU(\k):X_\a\in U_\a^\k\leftrightarrow\delta=\a .$$
\endproclaim

\demo{Proof} Let $X\in U_\a^\k$ be such that $X\notin U_\nu^\k$
for $\nu<\a .$ Put
$$f_\a(\xi)=\sup\{\eta\leq\oU(\xi);\forall\eta'<\eta : 
  X\cap\xi \notin U^\xi_{\eta'} \} .$$
Let $\delta\geq\a ,$ $j:V\ra M=\Ult(V,U_\delta^\k)$ then
$$M\models (jf_\a)(\k)=\sup\{\eta\leq o^{j\Cal U}(\k)=\delta;
  \forall \eta'<\eta : X\notin U^\k_{\eta'} \} .$$
Hence $(jf_\a)(\k)=\a$ since $X\in U_\a^\k$ and $X\notin U_{\eta'}^\k$
for $\eta'<\a .$ 
Put $Y=\{\xi<\k; f_\a(\xi)\geq\oU(\xi)\},$ then
$Y\in U_\delta^\k$ iff $\delta\leq \a .$
Consequently $X\cap Y$ separates $U_\a^\k$ from all the other
measures. \qed
\enddemo

It means that if there are no repeat points on $\k$ in $\Uarr$
then the measures on $\k$ are separable.

\proclaim{Lemma 3.8} If $\nu>\a$ has an $\a$-canonical function
and $\a$ is not a repeat point then $\nu$ is not a repeat point.
\endproclaim

\demo{Proof} Let $f_\a$ be the $\a$-canonical function for $\a$
defined above, and  $f_\nu$ an $\a$-canonical function for $\nu .$ 
Put $A=\{\xi<\k;\oU(\xi)>f_\a(\xi)\},$
$$B=A\cap \{\xi<\k; \oU(\xi)=f_\nu(\xi)\} .$$
Then $A\in U_\delta^\k$ iff $\delta>\a ,$ and so $B\in U_\nu$
separates $U_\nu$ from all the other measures. \qed
\enddemo 

Consequently if $\a$ is not a repeat point then $\a+1$ is not a repeat
point as $f_{\a+1}(\xi)=f_\a(\xi)+1$ is an $\a$-canonical function for
$\a+1 .$ The following is a joint result with J.Zapletal.

\proclaim {Proposition 3.9}
If $\a < \Theta$ then $j_{U^\k_\a}(\k) < \Theta .$
\endproclaim

\demo{Proof}
Put $\gamma = j_{U^\k_\a}(\k),$ $\a <\Theta .$
Let $f_\a$ be the $\a$-canonical function for $\a .$
We prove that for any $\nu <\gamma$, $U_\nu$ is not a repeat point.
Let $g_\nu : \k \rightarrow \k$ be a function such that
$(j_{U^\k_\a}g_\nu)(\k)=\nu .$
Define for $\lambda <\k$
$$
\tilde g_\nu(\lambda)= \cases (j_{U^\lambda_{f_\a(\lambda)}}(g_\nu
                               \restriction \lambda))(\lambda)
                               &\text{ if } o^{\Cal U}(\lambda)>f_\a(\lambda)\\
                               g_\nu(\lambda)
                               &\text{ otherwise .} \endcases
$$
I claim that $\tilde g_\nu$ is an $\a$-canonical function for $\nu .$
By the choice of $g_\nu$: $(j_{U_\a^\k}\tilde g_\nu)(\k)=\nu .$
Let $\delta \in (\a, o^{\Cal U}(\k))$ then
$$(j_{U^\k_\delta}\tilde g_\nu)(\k) = 
(j^{M_\delta}_{U^\k_\a}((j_{U^\k_\delta}g_\nu)
   \restriction \k)(\k) = (j_{U^\k_\a} g_\nu)(\k) = \nu $$
since $(j_{U^\k_\delta}g_\nu)\restriction \k=g_\nu$ and 
$j^{M_\delta}_{U^\k_\a} = j_{U^\k_\a} \restriction M_\delta ,$
where $M_\delta = \Ult(V,U^\k_\delta) .$
According to 3.8 $U^\k_\nu$ is not a repeat point.

Consequently $\gamma \leq \Theta .$ Note that in general 
$U \vartriangleleft W$ implies that $j_U(\k) < j_W(\k) .$
Since $\a+1 < \Theta$ we see that $\Theta \geq 
j_{U^\k_{\a+1}}(\k)>j_{U^\k_\a}(\k) .$
\qed
\enddemo

\proclaim {Remark} 
$\Theta$ is a limit ordinal between $\k^+$ and $\k^{++}$ of cofinality
$\k^+ .$ Moreover  if $F$ is a $\Sigma_1$ ordinal operation and
$\bar \a < \Theta$ then $F(\bar \a) < \Theta .$
\endproclaim

\demo{Proof}
Let $F$ be given  by a $\Sigma_1$ formula $\varphi (\bar \a, \beta)$
such that $ZFC\vdash$ ``$\varphi(\bar \a, \beta)$ is a function
of $\bar \a.$'' Let $\bar \a < \Theta$, find a $\delta < \Theta$
such that $\bar \a < j_{U_\delta}(\k)=\tilde \k$
($\delta > \bar \a$ is sufficient).
Since $M_\delta \models \text{``$\tilde \k$ is inaccessible,''}$
i.e.$V_{\tilde \k}^{M_\delta} \models ZFC ,$
it is necessary that
$$M_\delta \models \forall \bar \beta <\tilde \k :\;F^{V_{\tilde \k}}
                     (\bar \beta)<\tilde \k .$$
Moreover $M_\delta \models $``$ \forall\bar \beta<\tilde \k:
F^{V_{\tilde \k}}(\bar \beta) = F(\bar \beta)$'' and
$F^{M_\delta}=F^V$ as $F$ is a $\Sigma_1$ function.
Hence $F(\bar \a)<\tilde \k<\Theta.$
\qed
\enddemo

In particular $\Theta$ is inaccessible by primitive recursive
ordinal operations.

Finally we are ready to prove theorem 2.

\demo {Proof of theorem 2}
Given a  well-founded poset $P$ enumerate $P$ as in the proof of 3.6,
$P=\{p_\a; \a < \mu \} ,$ $\mu<\k^{++} .$
Then use 3.6 to find a closed system of measures
$\langle S_\lambda ;\lambda \leq \k \rangle$ such that
$P\cong (S_\k,\tS ),$ and a coherent sequence 
$\langle T_\lambda ; \lambda \leq \k \rangle$ such that $T_\k = S_\k$ and
the first measure $U_0$ satisfies $j_{U_0}(\k) > \mu .$
Then by  3.9 $T_\k$ cannot have a repeat point
as $o_T(\k)<j_{U_0}(\k) .$
Consequently the measures in $S_\k=T_\k$ are separable .
\qed
\enddemo

\enddocument
\end